\documentclass[12pt]{amsart}
\usepackage{graphicx}
\usepackage{amsmath}
\usepackage{amssymb}
\usepackage[numbers]{natbib}

\providecommand{\ph}[1]{\Phi_{#1}}
\providecommand{\pharg}[2]{\Phi_{#1}\left(#2\right)}
\providecommand{\stirlingfirst}[2]{\genfrac{[}{]}{0pt}{}{#1}{#2}}
\providecommand{\stirlingsecond}[2]{\genfrac{\{}{\}}{0pt}{}{#1}{#2}}

\theoremstyle{plain}
\newtheorem{thm}{Theorem}[section]
\newtheorem{lem}[thm]{Lemma}

\newtheorem{prop}[thm]{Proposition}

\theoremstyle{definition}
\newtheorem{defn}{Definition}
\newtheorem{exam}{Example}

\begin{document}
\title{On Separating Families of Bipartitions}

\author[T.\ Toda]{Takahisa Toda}
\address[T.\ Toda]{Graduate School of Human and Environmental Studies, Kyoto University, Kyoto, 606-8501 Japan.}
\email{toda.takahisa@gmail.com}

\author[I.\ Vigan]{Ivo Vigan}
\address[I.\ Vigan]{Department of Computer Science, City University of New York, The Graduate Center, New York, NY 10016, USA.}
\email{ivigan@gc.cuny.edu}
\thanks{The second author is supported by NSF grant 1017539.}

\keywords{set partition, Stirling number, shattering problem, spanning tree, hypergraph}
\subjclass[2010]{05A18}

\begin{abstract}
In this paper, we focus on families of bipartitions, i.e.\ set partitions consisting of at most two components.
We say that a family of bipartitions is a \emph{separating} family for a set $S$ if every two elements in $S$ can be separated by some bipartition.
Furthermore, we call a separating family \emph{minimal} if no proper subfamily is a separating family.
We characterize the set of all minimal separating families of maximum size for arbitrary set $S$ as the set of all spanning trees on $S$ and enumerate minimal separating families of maximum size.
Furthermore, we enumerate separating families of arbitrary size, which need not be minimal.
\end{abstract}

\maketitle

\section{Introduction}
A \emph{partition} $P$ of a set $S$ is a collection $\{U_{1},\ldots,U_{k}\}$ of pairwise disjoint nonempty subsets of $S$ with $S=\bigcup_{1\leq i\leq k}U_{i}$.
Each member $U_{i}$ is called a \emph{component} of $P$.
We mean by a \emph{bipartition} a set partition consisting of at most two components.
We call a bipartition \emph{proper} if it consists of exactly two components.

Any hypergraph $(V,E)$ can be considered as a family of bipartitions over $V$ because hyperedges $e$ in $E$ determine bipartitions of the form $\{e, V\setminus e\}$.
Bipartitions appear in multiple areas of computer science, one of them being the state-assignment problem in circuit design~\cite[chap.~12]{kohavi:automata}, where assignments can be considered as bipartitions.
Various separation algorithms for planar point sets by lines have been studied (see \cite{gruia:parallel-lines}, \cite{devillers:separation}, \cite{freimer:subdivisions}, \cite{freimer:shattering} and \cite{nandy:shattering}).
Toda~\cite{toda:partition} studied partitions of colored points by hyperplanes and proved a colorful Kirchberger-type Theorem.

In this paper we study \emph{separating} families of bipartitions as defined below.
We say that a partition $P$ \emph{cuts} two elements if they are in different components of $P$.
\begin{defn}
Let $S$ be a finite set.
A family $\mathcal{P}$ of bipartitions of $S$ is called a \emph{separating} family for $S$ if every two elements of $S$ can be cut by some bipartition in $\mathcal{P}$.
A separating family $\mathcal{P}$ for $S$ is \emph{minimal} if no proper subfamily of $\mathcal{P}$ is a separating family for $S$.
\end{defn}

\begin{exam}\label{exam:bipartitions}
Let $S=\{1,2,3,4\}$.
Let $P_{1},P_{2},Q_{1},Q_{2},Q_{3}$ be the bipartitions defined as:
\begin{align*}
P_{1} & = \{\{1, 2\}, \{3, 4\}\}, & Q_{1} & = \{\{1\}, \{2, 3, 4\}\}, \\
P_{2} & = \{\{1, 3\}, \{2, 4\}\}, & Q_{2} & = \{\{1, 2\}, \{3, 4\}\}, \\
& & Q_{3} & = \{\{1, 2, 3\}, \{4\}\}.
\end{align*}
Then the family of bipartitions $\{P_{1},P_{2}\}$ is a minimal separating family of minimum size for $S$, while $\{Q_{1},Q_{2},Q_{3}\}$ is that of maximum size.
\end{exam}

It is straightforward to check that the maximum size of a minimal separating family for an $n$-element set is $n-1$.

{\bf Our Contribution.}
We characterize the set of all minimal separating families of maximum size for arbitrary set $S$ as the set of all spanning trees on $S$ and enumerate minimal separating families of maximum size.
Furthermore, we calculate the number $\tau_{n,k}$ of separating families of $k$ arbitrary bipartitions over an $n$-element set and the number $\sigma_{n,k}$ of separating families of $k$ proper bipartitions over an $n$-element set.
Note that we distinguish between separating families of arbitrary bipartitions and those of proper bipartitions, both of which would be of interest to study.

This paper is organized as follows.
In Section~\ref{sect:bijection1} we enumerate all minimal separating families of maximum size for arbitrary set $S$ by obtaining a bijection onto the set of all spanning trees on $S$.
In Section~\ref{sect:arbitsize} we extend this analysis and calculate the number of separating families of arbitrary size.

\section{Minimal separating families of maximum size and spanning trees}\label{sect:bijection1}
In this section we enumerate all minimal separating families of maximum size for arbitrary set $S$ by obtaining a bijection onto the set of all spanning trees on $S$.

\begin{defn}
Let $S$ be a nonempty finite set.
For any minimal separating family $\mathcal{P}$ of bipartitions of $S$ we associate a graph denoted by $\pharg{S}{\mathcal{P}}$ such that the vertices of the graph correspond to the elements in $S$ and two vertices are adjacent if the corresponding elements in $S$ are cut by exactly one bipartition in $\mathcal{P}$.
\end{defn}

\begin{figure*}[ht]
  \begin{minipage}{0.45\hsize}
    \begin{center}
      \includegraphics[height=4cm]{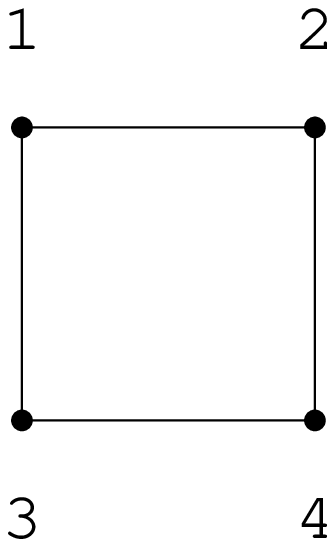}
    \end{center}
  \end{minipage}
  \begin{minipage}{0.45\hsize}
    \begin{center}
      \includegraphics[height=4cm]{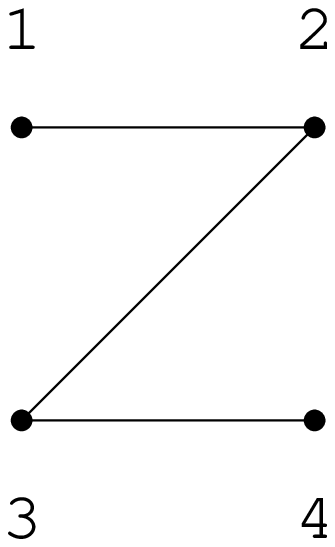}
    \end{center}
  \end{minipage}
  \caption{The graphs induced by the two minimal separating families of bipartitions $\{P_{1},P_{2}\}$ and $\{Q_{1},Q_{2},Q_{3}\}$ given in Example~\ref{exam:bipartitions}.}
\end{figure*}

\begin{thm}
Let $S$ be a nonempty finite set.
The mapping $\ph{S}\colon \mathcal{P}\mapsto\pharg{S}{\mathcal{P}}$ is a bijection from the set of all minimal separating families of maximum size for $S$ to the set of all spanning trees on $S$.
\end{thm}

\begin{proof}
Let $n$ denote the size of $S$, and remember that the maximum size of a minimal separating family for an $n$-element set is $n-1$.
We show that, for any minimal separating family $\mathcal{P}$ of maximum size for $S$, the associated graph $\pharg{S}{\mathcal{P}}$ is a spanning tree on $S$.
Let $G$ be a spanning subgraph of $\pharg{S}{\mathcal{P}}$ such that every bipartition in $\mathcal{P}$ has exactly one edge in $G$ whose end vertices are cut by it.
Clearly such a subgraph exists.
Assume that $G$ contains a cycle.
Let $C$ be such a cycle.
For any edge $e$ in $C$, the end vertices of $e$ can be connected by a path in $C$ which does not contain $e$.
By construction, in such a path, there are no edges whose end vertices can be cut by $P_{e}$.
This implies that the end vertices of $e$ are contained in the same component of $P_{e}$.
However, this contradicts that the end vertices of $e$ are cut by $P_{e}$.
Thus $G$ is a forest.
Since the vertex set of $G$ is of size $n$ and the edge set of $G$ is of size $n-1$, it follows that $G$ is a spanning tree on $S$.

Assume that $G$ is a proper subgraph of $\pharg{S}{\mathcal{P}}$.
For any edge $f$ in $\pharg{S}{\mathcal{P}}$ but not in $G$, the graph $G\cup\{f\}$ contains a cycle.
By a similar argument as above we can derive a contradiction.
Thus we have proved that $\pharg{S}{\mathcal{P}}=G$ and $\pharg{S}{P}$ is indeed a spanning tree on $S$.
It is straightforward to see that $\ph{S}$ is a one-to-one and onto mapping.
\end{proof}

Since Cayley's formula~\cite{cayley:trees} states that the number of spanning trees in the complete graph on $n$ labeled vertices is $n^{n-2}$, we immediately obtain the following enumeration result:

\begin{thm}
The number of minimal separating families of maximum size for an $n$-element set is $n^{n-2}$.
\end{thm}

\section{Enumerating separating families of arbitrary size}\label{sect:arbitsize}
In this section we study separating families of arbitrary size, which need not be minimal, and calculate the number of such separating families.

Throughout this section we assume without loss of generality that $S=\{1,\ldots,n\}$.

\begin{defn}
For any bipartition $P$ of $S$, we define $b(P)$ to be the vector of length $n$ whose $i$-th coordinate is given by
\begin{align*}
b_{i}(P) =\left\{
\begin{array}{ll}
1 & \text{if $P$ cuts $1$ and $i$}\\
0 & \text{otherwise.}
\end{array}
\right.
\end{align*}
Any $k$-tuple $P = (P_{1}, \ldots, P_{k})$ of bipartitions then corresponds to an $n \times k$ matrix whose $j$-th column vector is $b(P_{j})$.
We denote this matrix by $M_{P}$.
\end{defn}

We distinguish between a family of bipartitions and a tuple of bipartitions by denoting one as $\mathcal{P}$, $\mathcal{Q}$, \textit{etc} and the other as $P$, $Q$, \textit{etc}.

\begin{exam}
Let $S=\{1,2,3,4\}$.
Let $P=(P_{1},P_{2})$ and $Q=(Q_{1},Q_{2},Q_{3})$ be the two tuples whose bipartitions are given in Example~\ref{exam:bipartitions}.
Then we have
\begin{align*}
M_{P} & =
\begin{pmatrix} 
0 & 0 \\
0 & 1 \\
1 & 0 \\
1 & 1 \\
\end{pmatrix}, & 
M_{Q} & =
\begin{pmatrix} 
0 & 0 & 0 \\
1 & 0 & 0 \\
1 & 1 & 0 \\
1 & 1 & 1
\end{pmatrix}.
\end{align*}
\end{exam}

The following lemma is straightforward to verify and we omit the proof.
Note that arbitrary bipartitions of $S$ may include the partition $\{S\}$ consisting of the single component $S$.

\begin{lem}\label{lem:matrix}
The mapping $P\mapsto M_{P}$ is a bijection from the set of all $k$-tuples of arbitrary bipartitions of $S$ to the set of all $(0,1)$-matrices of size $n\times k$ such that the entries in the first row are all $0$.
\end{lem}

\begin{lem}\label{lem:distinct}
Let $\mathcal{P}$ be a family of bipartitions of $S$, and let $P$ be a tuple obtained from $\mathcal{P}$ by ordering its members.
Then $\mathcal{P}$ is a separating family for $S$ if and only if every two of the row vectors of $M_{P}$ are distinct.
\end{lem}

\begin{proof}
It is to see that for any two elements $i,j\in S$, there is a bipartition in $\mathcal{P}$ cutting them if and only if the $i$-th row vector and the $j$-th row vector of $M_{P}$ are distinct.
From this observation we immediately obtain the lemma.
\end{proof}

\begin{prop}\label{prop:minsize}
The minimum size of a separating family for an $n$-element set is $\lceil\log_{2}n\rceil$.
\end{prop}

\begin{proof}
It immediately follows from Lemma~\ref{lem:distinct} that the minimum size of a separating family is at least $\lceil \log_{2} n \rceil$, since otherwise at least two row vectors would coincide.
To see that the minimum size is at most $\lceil \log_{2} n\rceil$, we construct a tuple $P$ of bipartitions in such a way that the row vectors of the corresponding matrix $M_{P}$ are all distinct.
To achieve this, a row length of $\lceil \log_{2} n \rceil$ is sufficient and by Lemma \ref{lem:distinct}, the family consisting of the entries of $P$ is a separating family.
\end{proof}

On the other hand, the maximum size of a separating family for an $n$-element set is $2^{n-1}$ because a separating family of maximum size contains all possible bipartitions.

Let us denote by $\stirlingsecond{k}{i}$ the number of partitions of a $k$-element set into $i$ nonempty subsets.
This number is known as a \emph{Stirling number of the second kind} (see \cite[\S6.1]{graham:math}).

\begin{lem}\label{lem:stirlingsecond}
The number of sequences of length $k$ containing each element from a set of $i$ symbols at least once is $i!\,\stirlingsecond{k}{i}$.
\end{lem}

\begin{proof}
If we identify distinct positions $1,\ldots,k$ in a sequence of length $k$ with $k$ distinct objects, then there are $\stirlingsecond{k}{i}$ many ways to distribute $k$ distinct objects into $i$ indistinguishable boxes, with no box left empty.
If the boxes are distinguishable, then there are $i!\,\stirlingsecond{k}{i}$ many possibilities.
\end{proof}

Let us denote by $\tau_{n,i}$ the number of separating families of $i$ arbitrary bipartitions over an $n$-element set.
 
\begin{lem}\label{lem:arbitsize}
The following equation holds for any number $n,k$ with $2\leq n$ and $1\leq k\leq 2^{n-1}$:
\begin{align*}
\sum_{i=1}^{k}i!\,\stirlingsecond{k}{i}\,\tau_{n,i}=(2^{k}-1)(2^{k}-2)\cdots(2^{k}-n+1).
\end{align*}
\end{lem}

\begin{proof}
From Proposition~\ref{prop:minsize} we obtain that $\tau_{n,k}=0$ for $k<\lceil\log_{2}n\rceil$.
Thus we can assume that $k\geq\lceil\log_{2}n\rceil$.
There are $\binom{2^{k}-1}{n-1}\,(n-1)!$ many $(0,1)$-matrices of size $n \times k$ such that the entries in the first row are all $0$ and every two row vectors are distinct.
Among them, by Lemma~\ref{lem:matrix},~\ref{lem:distinct}, and \ref{lem:stirlingsecond}, there are $i!\,\stirlingsecond{k}{i}\,\tau_{n,i}$ many matrices such that the number of different column vectors is $i$.
This holds since each such matrix corresponds to a tuple of length $k$ containing each member from a separating family of size $i$ at least once.
Thus we obtain the equation stated in the lemma.
\end{proof}

Let us denote by $\stirlingfirst{k}{i}$ the number of permutations of $k$ elements which contain exactly $i$ permutation cycles.
This number is known as an \emph{unsigned Stirling number of the first kind} (see \cite[\S6.1]{graham:math}).

\begin{thm}\label{thm:tau1}
The number $\tau_{n,k}$ of separating families of $k$ arbitrary bipartitions over an $n$-element set with $2\leq n$ and $1\leq k\leq 2^{n-1}$ is
\begin{align*}
\tau_{n,k} = \frac{(n-1)!}{k!}\sum_{i=1}^{k}(-1)^{k-i}\,\stirlingfirst{k}{i}\,\binom{2^{i}-1}{n-1}.
\end{align*}
\end{thm}

\begin{proof}
The inversion formula for Stirling numbers (see \cite[\S3.1]{berge:comb}) states that if two integer sequences $\{a_{i}\}_{1\leq i\leq n_{0}}$ and $\{b_{i}\}_{1\leq i\leq n_{0}}$ satisfy $b_{k}=\sum_{i=1}^{k}\stirlingsecond{k}{i}\,a_{i}$, they also satisfy $a_{k}=\sum_{i=1}^{k}(-1)^{k-i}\,\stirlingfirst{k}{i}\,b_{i}$.
From Lemma~\ref{lem:arbitsize} we obtain the theorem.
\end{proof}

\begin{prop}
The number of separating families of minimum size for an $n$-element set is 
\begin{align*}
\frac{(n-1)!}{\lceil\log_{2}n\rceil!}\,\binom{2^{\lceil\log_{2}n\rceil}-1}{n-1}.
\end{align*}
\end{prop}

\begin{proof}
The number stated in the proposition can be calculated by the formula obtained in Theorem~\ref{thm:tau1}.
Note that if the index $i$ of the formula is in the range $1\leq i<\lceil\log_{2}n\rceil$, then $2^{i}<n$ and thus $\binom{2^{i}-1}{n-1}=0$.
\end{proof}

This number appears as the number of state assignments of an $n$-state machine in switching theory (see \cite[chap.~12]{kohavi:automata}).
See also \cite[\S16.2]{jukna:comb} for switching networks and DeMorgan formulas, \textit{etc}.

Let us denote by $\sigma_{n,k}$ the number of separating families of $k$ proper bipartitions over an $n$-element set.
We say that a bipartition of a set $S$ is \emph{trivial} if it has the form $\{S\}$.

\begin{lem}\label{lem:sigma1}
The equation $\sigma_{n,k}+\sigma_{n,k-1} = \tau_{n,k}$ holds for any numbers $n,k$ with $2\leq n$ and $2\leq k\leq 2^{n-1}$.
\end{lem}

\begin{proof}
Observe that any separating family of size $k$ for an $n$-element set $S$ which contains the trivial bipartition $\{S\}$ consists of $(k-1)$ proper bipartitions, which form a separating family for $S$, plus $\{S\}$.
Thus $\sigma_{n,k-1}$ coincides with the number of separating families of size $k$ for an $n$-element set which contain the trivial bipartition.
Thus we obtain the lemma.
\end{proof}

\begin{lem}\label{lem:stirlingfirst}
\begin{align*}
\stirlingfirst{k+1}{i+1} &= k!\,\sum_{j=i}^{k}\frac{1}{j!}\,\stirlingfirst{j}{i}
\end{align*}
\end{lem}

\begin{proof}
It holds (see \cite[\S6.1]{graham:math}) that $\stirlingfirst{k+1}{i+1} = k!\sum_{j=0}^{k}\stirlingfirst{j}{i}/j!$.
Since $\stirlingfirst{j}{i}=0$ for $0\leq j< i$, we obtain the equation.
\end{proof}

\begin{thm}\label{thm:sigma1}
The number $\sigma_{n,k}$ of separating families of $k$ proper bipartitions over an $n$-element set with $2\leq n$ and $1\leq k< 2^{n-1}$ is 
\begin{align*}
\sigma_{n,k} & = \frac{(n-1)!}{k!}\sum_{i=1}^{k}\,(-1)^{k-i}\,\stirlingfirst{k+1}{i+1}\,\binom{2^{i}-1}{n-1}.
\end{align*}
\end{thm}

\begin{proof}
Using Lemma~\ref{lem:sigma1} and~\ref{lem:stirlingfirst}, we can calculate $\sigma_{n,k}$ as follows:
\begin{align*}
\sigma_{n,k} & = \sum_{j=1}^{k}(-1)^{k-j}\tau_{n,j}\\
& = (n-1)!\sum_{j=1}^{k}\sum_{i=1}^{j}\frac{(-1)^{k-i}}{j!}\,\stirlingfirst{j}{i}\,\binom{2^{i}-1}{n-1}\\
& = (n-1)!\sum_{i=1}^{k}\sum_{j=i}^{k}\frac{(-1)^{k-i}}{j!}\,\stirlingfirst{j}{i}\,\binom{2^{i}-1}{n-1}\\
& = (n-1)!\sum_{i=1}^{k}\frac{(-1)^{k-i}}{k!}\,\stirlingfirst{k+1}{i+1}\,\binom{2^{i}-1}{n-1}.
\end{align*}
\end{proof}

For any $(0,1)$-matrix $M$ of size $n\times k$ whose entries in the first row and the first column are all zero, the transposed matrix $M^{t}$ represents a family of bipartitions of a $k$-element set which contains the trivial bipartition.

\begin{exam}
Let $S=\{1,2,3,4\}$.
Let $P=(\{S\},P_{1},P_{2})$ be the tuple whose bipartitions are given in Example~\ref{exam:bipartitions}.
Then we have
\begin{align*}
M_{P} & =
\begin{pmatrix} 
0& 0 & 0 \\
0& 0 & 1 \\
0& 1 & 0 \\
0& 1 & 1 \\
\end{pmatrix}, & 
M_{P}^{t} & =
\begin{pmatrix} 
0 & 0 & 0 & 0 \\
0 & 0 & 1 & 1 \\
0 & 1 & 0 & 1 \\
\end{pmatrix}.
\end{align*}
Thus, we obtain the following $4$ bipartitions over $\{1,2,3\}$:
\begin{align*}
\{\{1,2,3\}\},\ \{\{1,2\},\{3\}\},\ \{\{1,3\},\{2\}\},\ \{\{1\},\{2,3\}\}.
\end{align*}
\end{exam}

\begin{lem}\label{lem:transpose}
The equation $\sigma_{n,k-1}\,(k-1)! = \sigma_{k,n-1}\,(n-1)!$ holds for any numbers $n,k$ with $2\leq n$ and $2\leq k\leq 2^{n-1}$.
\end{lem}

\begin{proof}
The LHS of the equation is the number of $(0,1)$-matrices $M$ of size $n\times k$ such that the entries in the first row and the first column are all zero and $M$ corresponds to a separating family of size $k$ for an $n$-element set.
On the other hand, the RHS is the number of $(0,1)$-matrices $M'$ of size $k\times n$ such that the entries in the first row and the first column are all zero and $M'$ corresponds to a separating family of size $n$ for a $k$-element set.
The transpose operation induces a bijection between them, and thus we obtain the equation.
\end{proof}

From Theorem~\ref{thm:sigma1} and Lemma~\ref{lem:transpose}, one obtains the following theorem.

\begin{thm}\label{thm:sigma2}
The number $\sigma_{n,k}$ of separating families of $k$ proper bipartitions over an $n$-element set with $2\leq n$ and $1\leq k< 2^{n-1}$ is 
\begin{align*}
\sigma_{n,k} &=\sum_{i=1}^{n-1}(-1)^{n-1-i}\,\stirlingfirst{n}{i+1}\,\binom{2^{i}-1}{k}.
\end{align*}
\end{thm}

\begin{proof}
\begin{align*}
\sigma_{n,k} &= \frac{(n-1)!}{k!}\,\sigma_{k+1,n-1}\\
&=\frac{(n-1)!}{k!}\,\frac{k!}{(n-1)!}\,\sum_{i=1}^{n-1}(-1)^{n-1-i}\,\stirlingfirst{n}{i+1}\,\binom{2^{i}-1}{k}.
\end{align*}
\end{proof}

From Lemma~\ref{lem:sigma1} and Theorem~\ref{thm:sigma2}, one obtains the following theorem.

\begin{thm}\label{thm:tau2}
The number $\tau_{n,k}$ of separating families of $k$ arbitrary bipartitions over an $n$-element set with $2\leq n$ and $2\leq k< 2^{n-1}$ is 
\begin{align*}
\tau_{n,k}&= \sum_{i=1}^{n-1}(-1)^{n-1-i}\,\stirlingfirst{n}{i+1}\,\binom{2^{i}}{k}.
\end{align*}
\end{thm}

\begin{proof}
\begin{align*}
\tau_{n,k}&= \sigma_{n,k}+\sigma_{n,k-1}\\
&= \sum_{i=1}^{n-1}(-1)^{n-1-i}\,\stirlingfirst{n}{i+1}\,\left\{\binom{2^{i}-1}{k}+\binom{2^{i}-1}{k-1}\right\}\\
&= \sum_{i=1}^{n-1}(-1)^{n-1-i}\,\stirlingfirst{n}{i+1}\,\binom{(2^{i}-1)+1}{k}.
\end{align*}
\end{proof}

\begin{prop}
The minimum size of an element set for which there is a separating family of $k$ arbitrary bipartitions is $\lceil\log_{2}k\rceil+1$.
\end{prop}

\begin{proof}
For any number $k\ (\geq 1)$, let $n$ be any integer with $n<\lceil\log_{2}k\rceil+1$.
Then we obtain $k> 2^{n-1}$. 
Since the number of all bipartitions of an $n$-element set is $2^{n-1}$, there are no separating families of size $k$ for an $n$-element set.
Thus, the minimum size of an element set is at least $\lceil\log_{2}k\rceil+1$.
To see that the minimum size is at most $\lceil\log_{2}k\rceil+1$, let $n=\lceil\log_{2}k\rceil+1$.
Then we obtain $\lceil\log_{2}n\rceil\leq k\leq 2^{n-1}$.
Thus there is a separating family of size $k$ for an $n$-element set.
\end{proof}

\begin{prop}
The number of separating families of $k\ (\geq 2)$ arbitrary bipartitions whose element set is of minimum size is 
\begin{align*}
\binom{2^{\lceil\log_{2}k\rceil}}{k}.
\end{align*}
\end{prop}

\begin{proof}
The number stated in the proposition can be calculated by the formula obtained in Theorem~\ref{thm:tau2}.
Note that if the index $i$ of the formula is in the range $1\leq i<\lceil\log_{2}k\rceil$, then $2^{i}< k$ and thus $\binom{2^{i}}{k}=0$.
\end{proof}

The following two propositions can be proved in a similar way as above.

\begin{prop}
The minimum size of an element set for which there is a separating family of $k$ proper bipartitions is $\lceil\log_{2}(k+1)\rceil+1$.
\end{prop}

%\begin{proof}
%For any number $k\ (\geq 1)$, let $n$ be any integer with $n<\lceil\log_{2}(k+1)\rceil+1$.
%Then we obtain $k> 2^{n-1}-1$. 
%Since the number of all proper bipartitions of an $n$-element set is $2^{n-1}-1$, there are no separating families of $k$ proper bipartitions over an $n$-element set.
%Thus, the minimum size of an element set is at least $\lceil\log_{2}(k+1)\rceil+1$.
%To see that the minimum size is at most $\lceil\log_{2}(k+1)\rceil+1$, let $n=\lceil\log_{2}(k+1)\rceil+1$.
%Then we obtain $\lceil\log_{2}n\rceil\leq k\leq 2^{n-1}-1$.
%Thus there is a separating family of $k$ proper bipartitions over an $n$-element set.
%\end{proof}

\begin{prop}
The number of separating families of $k\ (\geq 1)$ proper bipartitions whose element set is of minimum size is 
\begin{align*}
\binom{2^{\lceil\log_{2}(k+1)\rceil}-1}{k}.
\end{align*}
\end{prop}

%\begin{proof}
%The number stated in the proposition can be calculated by the formula obtained in Theorem~\ref{thm:sigma2}.
%Note that if the index $i$ of the formula is in the range $1\leq i<\lceil\log_{2}(k+1)\rceil$, then $2^{i}-1< k$ and thus $\binom{2^{i}-1}{k}=0$.
%\end{proof}

\section*{Acknowledgments}
This work started when the authors joined the $14$th Korean Workshop on Computational Geometry 2011 and we would like to thank Professor Yoshio Okamoto for organizing it. Furthermore, we would like to thank Professor Peter Bra{\ss} for his advice and Dr.~Yasuhide Numata for some corrections.
Finally the first author is grateful to Professor Hideki Tsuiki for his encouragement.

\bibliographystyle{abbrvnat}
\bibliography{manuscript}
\end{document}